\newcommand{\Z}{{\mathbb Z}}
\newcommand{\Q}{{\mathbb Q}}
\newcommand{\R}{{\mathbb R}}

\newcommand{\ZH}{\Z[H_1(M)]}
\newcommand{\twid}{\widetilde}
\newcommand{\iso}{\approx}

\def\co{\colon\thinspace}
\documentclass[12pt]{amsart}
\usepackage{amssymb}
\usepackage{amsthm}
\usepackage{amsmath}
\usepackage{amscd}
\usepackage{graphicx}
\usepackage[all]{xy}
\usepackage[dvips,colorlinks=true]{hyperref}
\hypersetup{%
    pdftitle={Turaev Torsion and Cohomology Determinants for 3-Manifolds with Boundary},
    pdfauthor={Christopher Truman},
    pdfkeywords={pdf, latex, tex, ps2pdf, dvipdfm, pdflatex},
    bookmarksnumbered,
    pdfstartview={FitH},
    urlcolor=cyan,
}%

\newtheorem{theorem}{Theorem}[section]

\newtheorem{lemma}[theorem]{Lemma}

\DeclareMathOperator{\Det}{Det}
\DeclareMathOperator{\Tors}{Tors}
\DeclareMathOperator{\Hom}{Hom}

\DeclareMathOperator{\sign}{sign}

\DeclareMathOperator{\Tor}{Tor}

\begin{document}
\title[Turaev torsion and cohomology Determinants]{Turaev torsion and cohomology determinants for 3--manifolds with boundary}
\author{Christopher Truman}
\address{Mathematics Department\\University of Maryland\\College Park, MD  20742 USA}
\email{cbtruman@math.umd.edu}
\urladdr{http://www.math.umd.edu/\textasciitilde cbtruman/}
\begin{abstract}
We obtain generalizations of some results of Turaev from \cite{Tur:3d}.  Turaev's results relate leading order terms of the Turaev torsion of closed, oriented, connected $3$--manifolds
to certain ``determinants'' derived from cohomology operations such as the alternate trilinear form on the first
cohomology group given by cup product.  These determinants unfortunately do not generalize directly to
compact, connected, oriented $3$--manifolds with nonempty boundary, because one must incorporate the cohomology of the manifold relative to
its boundary.  We define the new determinants that will be needed,
and show that with these determinants enjoy a similar relationship to the one
given in \cite{Tur:3d} between torsion and the known determinants.  These definitions and results are given for integral cohomology, cohomology with coefficients in $\mathbb{Z}/r\mathbb{Z}$ for certain integers $r$, and for integral Massey products.
\end{abstract}
\thanks{The author would like to thank James Schafer and Vladimir Turaev for valuable input and inspiration}
\date{\today}
\maketitle
\begin{section}{Introduction}
\label{section-Introduction}
The Turaev torsion is a way of combining all of the ``interesting'' Reidemeister torsions into a single invariant.  Dimension $3$ is a particularly nice place to look at Reidemeister torsions, because the cellular complexes one obtains from triangulating a manifold are reasonably easy to understand (we will use the word ``manifold'' to mean a compact, connected, orientable, and smooth manifold).  In addition, Turaev torsion is related to many interesting invariants in dimension $3$, including the Seiberg-Witten invariant and the Casson-Walker-Lescop invariant.  For excellent discussions of these ideas, see \cite{Tur:3d}, \cite{Tur:Intro}, and \cite{Nicol-Reid3} (though one should note that Turaev refers to what we call ``Turaev torsion'' as ``maximal abelian torsion,'' and Nicolaescu refers to it as ``Reidemeister-Turaev torsion'').  A relationship between Turaev torsion and cohomology for closed $3$--manifolds is known (see \cite{Tur:3d} Chapter~III), and below we give the analogue for $3$--manifolds with boundary.

We will give a definition of Turaev torsion based upon the definition of Reidemeister torsion.  First, a brief review of Reidemeister torsion.  Recall the (commutative) Reidemeister torsion of a finite complex $X$ and a ring map $\varphi:\Z[H_1(X)]\rightarrow F$, where $F$ is a field, is an element $\tau^\varphi(X)\in F/(\pm \varphi(\Z[H_1(X)])$.  Given an Euler structure $e$ and a homology orientation $\omega$ (see \cite{Tur:3d} for definitions), one can define the refined $\varphi$--torsion $\tau^\varphi(X,e,\omega)\in F$, which is a unit if the $\varphi$--twisted complex $C^\varphi_*(X)$ is acyclic, and is equal to $0\in F$ if not.
\begin{subsection}{The Turaev Torsion}
\label{subsection-TurTor}
If $X$ is a finite connected CW--complex, the quotient ring
(i.e. the ring obtained by localizing at the multiplicative set of
non-zerodivisors) of the integral group ring of a finitely generated abelian
group splits as a direct
sum of fields (see \cite{Tur:3d} I.3.1 for details).  This isomorphism provides ring homomorphisms from $\Z[H_1(X)]$
to various fields.  Specifically, if we denote by $Q(H)$ the quotient ring of
$\Z[H]$ where $H=H_1(X)$, we have
the inclusion $\Z[H] \hookrightarrow Q(H)$.  There is an isomorphism $\Phi\co Q(H)
\stackrel{\iso}{\rightarrow} \bigoplus\limits_i F_i$ where each $F_i$
is a field and $i$ ranges
over a finite index set.  This isomorphism is defined, for example, in \cite{Tur:3d}, and is unique up to unique isomorphism (which, after re-ordering if necessary, will decompose along the direct sum as a component-wise isomorphism $F_i\rightarrow F_i'$) making the following diagram commute:
\[
\xymatrix{
 & \bigoplus\limits_i F_i\ar[dd]\\
Q(H)\ar[ur]\ar[dr] & \\
& \bigoplus\limits_i F_i'.\\
}
\] 
Then denote by $\varphi_i$ the map
$\Z[H]\rightarrow F_i$ consisting of the inclusion to $Q(H)$ followed
by the natural projection to $F_i$.  Then for any homology orientation
$\omega$ and Euler structure $e$, we
define the Turaev torsion $\tau(X,e,\omega)$ by
\[ \tau(X,e,\omega) = \Phi^{-1}\left ( \bigoplus_i
\tau^{\varphi_i}(X,e,\omega) \right ) \in Q(H).\] 

This definition does not depend on $\Phi$ (by the uniqueness of $\Phi$).  Henceforth, the symbol $\tau(X,e,\omega)$ will refer to the Turaev torsion of
$(X,e,\omega)$ unless otherwise specified.  The symbol
$\tau^\varphi(X,e,\omega)$ will still refer to $\varphi$--torsion.  

The set of Euler structures possesses a free, transitive $H_1(X)$--action, and we will use $-\omega$ to represent the opposite homology orientation to $\omega$.  One may easily show that
\begin{align*}
\tau(X,he,\omega)&=h\tau(X,e,\omega)\\ 
\tau(X,e,-\omega)&=-\tau(X,e,\omega).
\end{align*}
\end{subsection}
\begin{subsection}{Turaev Torsion and Alexander Invariants}
\label{section-Tors-vs-Alex}
This section repeats a result from \cite{Tur:3d} II.1.7.  Though this result is unproven in \cite{Tur:3d}, it is a simple exercise to prove it via the method of proof of \cite{Tur:3d} Theorem~II.1.2.  We will use the result in the proof of Theorem~\ref{theorem-IntCohom-vs-tors}.

Henceforth, we will often
need to strike a column from a matrix.  We will use the notation $A(r)$ for the
matrix obtained by striking the $r^\text{th}$ column from a matrix $A$.

\begin{theorem}[Turaev]
\label{theorem-tors-vs-Alex}
Suppose $M$ is a $3$--manifold with $\partial M \neq \varnothing$, $\chi(M)=0$.  Then there is an Euler structure $e$ such that for any homology orientation 
$\omega$, and for any $1 \leq r \leq m$,
\begin{equation}
\tau(M,e,\omega)(h_r-1)=(-1)^{m+r}\tau_0 \det (\Delta(r)).
\label{equation:3MbdTorCalc}
\end{equation}
\end{theorem}
\end{subsection}
\end{section}

\begin{section}{The Integral Cohomology Ring}
\label{section-IntCohom}
%
\begin{subsection}{Determinants}
\label{subsection-IntCohom-Dets}
Let $M$ be a 
$3$--manifold with boundary $\partial M \neq \varnothing$, and suppose
$\chi(M) = 0$.  Also assume $b_1(M) \geq 2$ so that $H^1(M,\partial M) \neq 0$. Let $f_M$ denote the map $H^1(M,\partial M)\times H^1(M)\times H^1(M)
\longrightarrow \Z$, defined by $(a,b,c)\mapsto \langle a \cup b \cup
c , [M] \rangle$, where $[M]$ is the fundamental class in
$H_3(M,\partial M)$ determined by the orientation and $\langle \cdot,\cdot\rangle$ denotes evaluation pairing.  This 
is alternate in the last two variables; $\langle a \cup b \cup c , [M] \rangle =
-\langle a \cup c \cup b , [M] \rangle$.  Our assumption $\chi(M) =
0$ implies that $H^1(M,\partial M)$ will have rank $b_1(M)-1$.  There is a notion of a ``determinant'' of an alternate trilinear form (see \cite{Tur:3d},
chapter III for a discussion of the determinant of the obvious
analogue of the above form when $M$ is closed), but because of the
difference in rank, we must have a new concept of determinant for a 
mapping such as the one above.  The determinant of an alternate
trilinear form on a free $R$--module is independent of basis up to
squares of units of $R$, so if $R=\Z$ it is independent of basis.
This will not be true of our determinant; however we will present a
sign-refined version based on a choice of homology orientation.  For
our usage, this is not more of a choice than we would normally make; if
we want sign-refined torsion, then we have already chosen a homology
orientation, and if we do not care about the sign of the torsion, we
can ignore the sign here as well. 

Let $R$ be a commutative ring with unit, and let
$K,L$ be finitely generated free $R$ modules of rank $n$ and $n-1$
respectively, where $n \geq 2$. Let $S =
S(K^*)$ be the symmetric algebra on $K^*$, where $K^* = \Hom_R(K,R)$.  Note if $\{a_i\}_{i=1}^{n}$ is
a basis of $K$ with dual basis $\{a_i^*\}_{i=1}^{n}$ of $K^*$, then
$S=R[a_1^*,\dots,a_{n}^*]$, the polynomial ring on
$a_1^*,\dots,a_{n}^*$, and the grading of $S$ corresponds to the usual
grading of a polynomial ring. 

Now let
$\{a_i\}_{i=1}^n,\{b_j\}_{j=1}^{n-1}$ be bases for $K,L$ respectively,
and let $\{a_i^*\}$ 
be the basis of $K^*$ dual to the basis $\{a_i\}$.  Let $f\co L\times
K\times K \longrightarrow R$ be an $R$--module homomorphism which is
skew-symmetric in the two copies of $K$; i.e.  for
all $y,z \in K, x\in L, f(x,y,z)=-f(x,z,y)$.  Let $g$ denote the
associated homomorphism $L\times K\longrightarrow K^*$ given by
$(g(x,y))(z) = f(x,y,z)$.  Next we state a Lemma defining the {\it determinant} of $f$ ($d$ in the Lemma), but first we set some notation: $[a'/a] \in R^\times$ is used to denote the
determinant of the change of basis 
matrix from $a$ to $a'$, and for a
matrix $A$, we will let $A(i)$ denote the matrix obtained by striking
out the $i^{\text{th}}$ column as above.
\begin{lemma}
\label{lemma-IntCohomDet}
Let $\theta$ denote the $(n-1 \times n)$ matrix over $S$ whose
$i,j^{\text{th}}$ entry $\theta_{i,j},$ is obtained by
$\theta_{i,j} = g(b_i,a_j)$.  Then there is a unique $d = d(f,a,b) \in
S^{n-2}$ such that for any $1 \leq i \leq n$,
\begin{equation}
\det\ \theta(i) = (-1)^{i}a_i^*d.
\label{equation-IntCohomDetDef}
\end{equation}
For any other bases $a',b'$ of $K,L$ respectively, we have 
\begin{equation}
d(f,a',b') = [a'/a][b'/b]d(f,a,b).
\label{equation-IntCohomCob}
\end{equation}
\end{lemma}
\begin{proof}
Let $\beta$ denote the $(n-1\times n)$ matrix with $\beta_{i,j} =
g(b_i,a_j)a_j^*$.  The sum of the columns of $\beta$ is zero; indeed, for any
$i$, the $i^{\text{th}}$ entry (of the column vector obtained by summing the columns of $\beta$) is given by:
\begin{equation*}
\sum\limits_{j=1}^n \beta_{i,j} = \sum\limits_{j=1}^n g(b_i,a_j)a_j^* =
\sum\limits_{j,k=1}^n f(b_i,a_j,a_k)a_j^*a_k^* = 0.
\end{equation*}
The last equality follows since the $f$ term is anti-symmetric in
$j,k$ and the $a$ term is symmetric.  It is then a simple algebraic fact that since $\beta$ is an $n-1\times n$ matrix whose columns sum to zero,
$(-1)^{i}\det\beta(i)$ is independent of $i$.  

Let $t_i = \det
\theta(i) \in S^{n-1}$.  It is clear that
\begin{equation*}
\det (\beta(i)) = t_i \prod_{k \neq i}a_k^*,
\end{equation*}
thus for any $i,p \leq n$, we have
\begin{eqnarray*}
(-1)^{i}t_ia_p^*\prod_{k=1}^n a_k^* & = & (-1)^{i}\det \beta(i)
a_p^*a_i^* \\
& = & (-1)^p\det \beta(p) a_i^*a_p^*\\ 
& = & (-1)^pt_pa_i^*\prod_{k=1}^n a_k^*.
\end{eqnarray*}
Since the annihilators of $a_k^*$ in $S$ are zero, we must have 
\begin{equation*}
(-1)^{i}t_ia_p^* = (-1)^pt_pa_i^*.
\end{equation*}
This means that $a_i^*$ divides $t_ia_p^*$ for all $p$, hence $a_i^*$
divides $t_i$.  Define $s_i$ by $t_i = s_ia_i^*$.  Note 
\begin{equation*}
(-1)^{i}s_ia_i^*a_p^* = (-1)^{i}t_ia_p^* =
(-1)^pt_pa_i^*=(-1)^ps_pa_p^*a_i^*,
\end{equation*}
hence $(-1)^{i}s_i$ is independent of $i$.  Let $d =
(-1)^{i}s_i$.  By definition, 
\begin{equation*}
(-1)^{i}\det \theta(i) = (-1)^{i}t_i = (-1)^{i}s_ia_i^* = a_i^*d.
\end{equation*}
This proves \eqref{equation-IntCohomDetDef}.

Now to prove the change of basis formula \eqref{equation-IntCohomCob}, we will first
show $d(f,a',b) = [a'/a]d(f,a,b)$.  Let $S_i$ be the $(n\times n-1)$ matrix obtained by
inserting a row of zeroes into the $(n-1\times n-1)$ identity matrix
as the $i^\text{th}$ row.  Then one may easily see for any $(n-1\times
n)$ matrix $A$, the matrix $A(i)$ (obtained by striking out the
$i^\text{th}$ column) can also be obtained as $A(i)=AS_i$.  Let
$S_i^+$ denote the $(n\times n)$ matrix obtained by appending a column
vector with a $1$ in the $i^\text{th}$ entry and zeroes otherwise on
to the right of $S_i$ (i.e. $S_i^+$ is the identity matrix with the $i^\text{th}$ column moved all the way to the right), and let $A_+^i$ denote the $(n\times n)$ matrix
obtained by appending a row vector with a $1$ in the $i^\text{th}$
entry and zeroes otherwise on to the bottom of $A$.  Note
\begin{align*}
\det(S_i^+) &= (-1)^{n+i}\\
\det(A_+^i) &= (-1)^{n+i}\det(A(i))\\
&\text{hence}\\
\det(A_+^iS_i^+) &= \det(AS_i) = \det(A(i)).
\end{align*}
Now let $(a'/a)$
denote the usual change of basis matrix so that $a'_i =
\sum\limits_{j=1}^n (a'/a)_{i,j} a_j$.   Now $\theta_{i,j} =
g(b_i,a_j)$, so let 
\[
\theta'_{i,j} = g(b_i,a'_j).
\]
One can easily show $\theta' = \theta\cdot (a'/a)^\mathsf{T}$, where $\mathsf{T}$ denotes transpose.  

Now we compute
\begin{align*}
\det\left(\theta_+^i (a'/a)^\mathsf{T} S_i^+\right) &=
\det\left(\theta_+^i\right)\det\left((a'/a)^\mathsf{T}\right)\det\left(S_i^+\right) \\
&=\det\left((a'/a)^\mathsf{T}\right) \det\left(\theta_+^i\right)\det\left(S_i^+\right) \\
&=\det\left(a'/a\right) \det\left(\theta_+^i\cdot S_i^+\right) \\
&= [a'/a]\det(\theta(i)) \\
&=[a'/a](-1)^i a_i^* d(f,a,b).
\end{align*}
To complete the proof, we multiply the argument of the determinant in the left-hand side of the above equation out.  Let
$e_i$ denote the row vector with a $1$ in the $i^\text{th}$ position and zeroes otherwise, i.e. the $i^\text{th}$ basis vector of $a$ as expressed in the $a$--basis, and let $r_i$ denote
the $i^\text{th}$ row of $(a'/a)^\mathsf{T}$ and $c_i$ denote the
$i^\text{th}$ column.  Then {\small
\begin{align*}
\det(\theta_+^i (a'/a)^\mathsf{T} S_i^+) &= \det\left[\left(\begin{smallmatrix}
\theta \\ e_i\end{smallmatrix}\right) (a'/a)^\mathsf{T} \left(\begin{smallmatrix}
S_i & e_i^\mathsf{T} \end{smallmatrix}\right)\right] \\
&= \det\left[\left(\begin{smallmatrix}
\theta \\ e_i\end{smallmatrix}\right) \left(\begin{smallmatrix}
(a'/a)^\mathsf{T}(i) & c_i\end{smallmatrix}\right)\right] \\
&= \det\left(\begin{smallmatrix} \theta (a'/a)^\mathsf{T}(i) & \theta c_i
\\ (a'/a)^\mathsf{T}(i)_i & (a'/a)^\mathsf{T}_{i,i} \end{smallmatrix}\right) \\
&= (-1)^{n-i}\det\left(\begin{smallmatrix}\theta (a'/a)^\mathsf{T} \\
r_i\end{smallmatrix}\right) \\
&= (-1)^{n-i}\det\left(\begin{smallmatrix}\theta' \\
r_i\end{smallmatrix}\right) \\
&= (-1)^{n-i}\sum\limits_{k=1}^n
(-1)^{n+k}(a'/a)^\mathsf{T}_{i,k}\det(\theta'(k)) \\
&= (-1)^{n-i}\sum\limits_{k=1}^n (-1)^{n+k}(a'/a)^\mathsf{T}_{i,k} (-1)^k
(a'_k)^* d(f,a',b) \\
&= (-1)^{n-i}\sum\limits_{k=1}^n (-1)^{n+k}(a^*/(a')^*)_{i,k} (-1)^k
(a'_k)^* d(f,a',b) \\
&= (-1)^i d(f,a',b) a_i^*.
\end{align*}}
So $d(f,a',b)-[a'/a]d(f,a,b)$ annihilates $a_i^*$ for each $i$, hence
is zero.

The computation for a $b$ change of basis is easier.  
Let $b'$ be
another basis for $L$ and let $(b'/b)$ denote the $b$ to $b'$ change
of basis matrix.  Let $\theta'$ denote the matrix $g(b_i',a_j)$, then
$\theta' = (b'/b)\theta$.  So $\theta' S_i = (b'/b)\theta S_i$, hence
$\det(\theta'(i)) = [b'/b]\det(\theta(i))$.
This proves $d(f,a,b') = [b'/b]d(f,a,b)$, and completes the proof of
\eqref{equation-IntCohomCob}. 
\end{proof}

\begin{subsubsection}{The Sign Refined Determinant}
In the case, $R=\Z$,  our determinant depends on the basis only by
its sign.  In this case, we can refine the determinant by a choice of
orientation of the $\R$--vector space $(K\oplus L)\otimes \R$.  Let
$\omega$ be such a choice of orientation.  Then define $\Det_\omega(f)
= \det(f,a,b)$ where $a,b$ are bases of $K,L$ respectively such that
the induced basis of $(K\oplus L)\otimes \R$ given by
extension of scalars
is positively oriented with respect to 
$\omega$.  Then $\Det_\omega(f)$ is well defined, and for any bases
$a',b'$, we have $\det(f,a',b') = \pm \Det_\omega(f)$ where the $\pm$
is chosen depending on whether $a',b'$ induces a positively or
negatively oriented basis of $(K\oplus L)\otimes \R$ with respect to
$\omega$.

Note that for $K=H^1(M)$, $L=H^1(M,\partial M)$ where $M$
is a compact connected oriented $3$--manifold with non-void boundary, a choice
of homology orientation will determine an orientation for $(K\oplus
L)\otimes \R$.  One simply says that $a,b$ is a positively oriented basis for $H^1(M,\R)\oplus H^1(M,\partial M;\R)$ with respect to a homology orientation $\omega$ if and only if $\{[\text{pt}] , a_1^*,\dots,a_n^*,b_1\cap [M],\dots,b_{n-1}\cap [M]\}$ is a positively oriented basis in homology (where $a^*$ is the basis of $H_1(M,\R)$ dual to $a$, $[M]$ is the fundamental class determined by the orientation, and $[\text{pt}]$ is the homology class of a point).

This sign refinement is identical to refining $\Det$ by the paired volume form associated to $\omega$, as defined below in \eqref{equation-volFormRefinedDet}.  Also, note that this sign refined determinant only depends on the homology orientation, not the orientation of $M$.
\end{subsubsection}
\end{subsection}
%
\begin{subsection}{Relationship to Torsion}
\label{subsection-IntCohom-Relns}
We use the above to relate the torsion to the
cohomology ring structure.  Let $T = \Tors(H_1(M))$ denote the torsion
subgroup of $H_1(M)$.  Note that this is isomorphic to the torsion
subgroup of $H_1(M,\partial M)$, so we will also denote the torsion subgroup of
$H_1(M,\partial M)$ by $T$.  Let $G=H_1(M)/T$, let $S(G)$ denote the
graded symmetric
algebra on $G$ and let $I$ denote the augmentation ideal in $\ZH$.
The filtration of $\ZH$ by powers of $I$ determines an associated
graded algebra $A = \bigoplus\limits_{\ell \geq 0} I^\ell/I^{\ell +
1}$.  Then there
is an additive homomorphism $q\co S(G) \longrightarrow A$
defined in \cite{Tur:3d}.  We repeat the definition here:
The map $h \mapsto h-1 \bmod I^2$ defines an additive homomorphism
$H_1(M) \longrightarrow I/I^2$.  This extends to a grading-preserving
algebra homomorphism $\twid{q}_{H_1(M)}\co S(H_1(M))\longrightarrow A$.  Any
section $s\co G\longrightarrow H_1(M)$ of the natural projection
$H_1(M)\longrightarrow G$
induces an algebra homomorphism $\twid{s}\co S(G)\longrightarrow S(H_1(M))$; set 
\begin{equation*}
q =
|T|\twid{q}_{H_1(M)}\twid{s}\co S(G)\longrightarrow A.
\end{equation*}
Then $q$ is grading preserving and is a $\Z$--module
homomorphism, and obviously satisfies the multiplicative formula
\begin{equation*}
q(a)q(b) = |T|q(ab).
\end{equation*}
The map $q$ does not depend on the choice of section $s$ (see \cite{Tur:3d}).

We are now ready to state the main result of this section:
\begin{theorem}
\label{theorem-IntCohom-vs-tors}
Let $f_M\co H^1(M,\partial M)\times H^1(M)\times H^1(M)\longrightarrow \Z$ be
the $\Z$--module homomorphism defined by
\begin{equation*}
f_M(x,y,z)=\langle x\cup y\cup z , [M]\rangle.
\end{equation*}
Let $n = b_1(M) \geq 2$,
let $I$ be
the augmentation ideal of $\ZH$, and let $e$ be any choice of Euler
structure on $M$ and $\omega$ be a homology orientation of $M$.
Then $\tau(M,e,\omega)\in I^{n-2}$
and:
\begin{equation}
\tau(M,e,\omega) \bmod I^{n-1} = q(\Det_{\omega} (f_M)) \in
I^{n-2}/I^{n-1}.
\label{equation-IntCohomTorsThm}
\end{equation}
\end{theorem}
That $\tau(M,e,\omega)\in I^{n-2}$ is proved in
\cite{Tur:3d}, Chapter II, the theorem is concerned with its image
modulo $I^{n-1}$; this is the ``leading term'' of the torsion in the
associated graded algebra $A$.  This proof is generally the method of \cite{Tur:3d} Theorem~2.2, though some adjustments must be made to incorporate the relative homology.
\begin{proof}
The first step is to arrange a handle decomposition coming from a $C^1$ triangulation to be in a convenient form for comparing the torsion to the cohomology. 

First, we arrange our decomposition for $M$ so that we have $(0)$ $3$--handles, $(m-1)$
$2$--handles, $(m)$ $1$--handles, $(1)$ $0$--handle, and this is Poincar\'{e} dual to a relative handle
decomposition for $(M,\partial M)$ with $(0)$ $0$--handles, $(m-1)$
$1$--handles, $(m)$ $2$--handles, $(1)$ $3$--handle.  
With these decompositions, we have the following
cellular chain complexes:
\begin{equation*}\begin{CD}
C_*(M): @. 0 @>>> \Z^{m-1} @>>> \Z^m @>0>> \Z\\
@. @VV\approx V @VV\approx V @VV\approx V @VV\approx V\\
C_*(M,\partial M):\ \ \ \  @. 0 @<<< \Z^{m-1} @<<< \Z^m @<<0< \Z.
\end{CD}\end{equation*}

We will refer to the handles as ``honest''
handles and ``relative'' handles; honest handles being from the
decomposition of $M$ and relative ones from the relative decomposition
of $(M,\partial M)$.  Later, we will explicitly give the $(m-1 \times m)$ matrix for $\partial_2$ of the honest decomposition (i.e. the only nonzero boundary in the honest complex).  

The core $0$--cell of the honest $0$--handle (of $M$) is a point, $u$, which
we will say is positively oriented.  At the same time we orient the
relative $3$--handle (of $(M,\partial M)$) with the positive orientation given by
the orientation of $M$.  Extend the core $1$--disks of the
honest $1$--handles to obtain loops in $M$ based at $u$, representing
$x_1,\dots,x_m \in \pi_1(M,u)$.  Since sliding the
$i^\text{th}$ honest $1$--handle over the $j^\text{th}$ honest $1$--handle
replaces $x_i$ with $x_ix_j$, and reversing orientation of the core
$1$--disk changes replaces $x_i$ with $x_i^{-1}$, we may perform handle moves to assume that the images of the 
homology classes of the first $n$ of the $x_i$'s form a basis of $G = H_1(M)/T$ and the rest of the classes end up in $T$.
For 
$i=1,\dots,m$, let $h_i$ be the image of $x_i$ in $H_1(M)$ under the Hurewizc map,
and $\twid{h_i} = h_i \bmod T$.  Thus
$\twid{h_1},\dots,\twid{h_n}$ is a basis of $G$ and $\twid{h_i} = 1$
for $i > n$.  Denote the 
dual basis of $H^1(M)$ by $h^*_1,\dots,h^*_n$, by definition, $\langle
h^*_i,\twid{h_j} \rangle = \delta_{i,j}$, where $\langle \cdot , \cdot
\rangle$ is evaluation pairing.

We now arrange the relative $1$--handles in a similar way;
so that the images of the first $n-1$ of them form a basis of $H_1(M,\partial
M)/T$ and the other $m-n$ of them end up in the
torsion subgroup.  Let $c$ denote the number of components of $\partial
M$, we will now describe a method to geometrically realize the splitting of $H_1(M,\partial M)$ as the image of $H_1(M)$ direct sum with $\Z^{c-1}$ generated by paths connecting boundary components.  We will arrange for the first
$c-1$ of the relative handles to connect a given boundary component to the other components, and the rest of the handles to represent loops based at a point on that component.  

Since $M$ is connected, the relative $1$--skeleton is path connected, which means boundary component has at least one relative $1$--handle with only one endpoint on that component.  Choose a base point in one boundary component, call that component $(\partial M)_0$.  By sliding relative handles over one another if necessary, we can arrange so that there is a relative $1$--handle connecting $(\partial M)_0$ to each other boundary component.  This is reasonably simple to do; given another boundary component, there is a path through the relative $1$--skeleton connecting that component to $(\partial M)_0$.  The first part of this path is a relative $1$--handle with one endpoint in the boundary component we would like to connect to $(\partial M)_0$, so one simply slides the other end of the relative $1$--handle along the path to $(\partial M)_0$, which is simply sliding it over other boundary components and other relative $1$--handles, until it connects our boundary component to $(\partial M)_0$.  Once we have fixed $c-1$ relative handles doing this (one for each boundary component which is not $(\partial M)_0$), slide all of the other relative handles along the fixed $c-1$ handles to obtain relative handles whose cores generate the image of the fundamental group of $M$ based at the chosen base point.  

To summarize, by handle moves, we arrange so that the first $c-1$ relative handles are paths connecting a chosen boundary component to all of the other boundary components, and the other $m-c$ relative handles have homology classes generating the image of $H_1(M)$ in $H_1(M,\partial M)$.  With these last handles, we 
may proceed as before in the discussion of honest handles; arrange so that the first $n-c$ of these handles will give us the
remaining free generators of $H_1(M,\partial
M)/T$ and the rest of them simply end up in $T$ (again by sliding
handles, since the only handles that we now need to slide represent loops
all based at the same point).  We will use similar notation, $k_i$
will denote the 
homology class of the $i^\text{th}$ handle and $\twid{k_i}=k_i \bmod
T$.  We
will denote the dual basis of $H^1(M,\partial M)$ by
$k^*_1,\dots,k^*_{n-1}$.  As before, the $\twid{k_i}$'s for $i \leq
n-1$ are generators of $H_1(M,\partial
M)/T$ and for $i > n-1,
\twid{k}_i = 1$.  Also, as before, $\langle k^*_i,\twid{k_j} \rangle =
\delta_{i,j}$. 

The attaching maps for the honest $2$--handles determine (up to
conjugation) certain elements $r_1,\dots,r_{m-1}$ of the free group
$F$ generated by $x_1,\dots,x_m$.  We now have $\pi_1(M)$ presented by
the generators $x_1,\dots,x_m$ and the relations $r_1,\dots,r_{m-1}$.

Now the cellular chain complex for $M$ is in a particularly convenient form for our purposes.  As usual, we use the notation $\partial_p$ to denote the boundary
map from dimension $p$ to dimension $p-1$.  Clearly $\partial_1$ is given by the zero map.  Let us denote the matrix of
$\partial_2$ by $(v_{i,j})$ where $1 \leq i \leq
m-1$ and $1 \leq j \leq m$.

Now for $1 \leq i \leq n-1$, the core $2$--disk of the $i^{\text{th}}$
honest $2$--handle represents a cycle in $C_2(M)$ (in fact, its homology class is $k_i^*\cap [M]$).  So it has boundary equal to zero, hence $v_{i,j} =
0$ for $i \leq n-1$ and all $j$.  We apply the same argument to the relative handles as follows:  the $j^\text{th}$ relative $2$--handle represents a homology class Poincar\'{e} dual to $h_i^*\cap [M]$, hence has boundary equal to zero, and $v_{i,j} =
0$ for all $i$ and $j\leq n$.  The result is that $v_{i,j} = 0$ except for the
bottom right hand $(m-n \times m-n)$ corner of the matrix; call this matrix $v$.  This
tells us that $\partial_2$ in the
complex for $M$ is given by
$\left(\begin{smallmatrix}0&0\\0&v\end{smallmatrix}\right)$. This $v$ is a square
presentation matrix for the torsion group $T$, thus $\det (v) = \pm |T|$.  Furthermore,  $r_1,\dots,r_{n-1} \in [F,F]$ since the first $n-1$ honest $2$--cells are cycles.  Essentially, sliding handles and relative handles amount to doing integral row and column operations to reduce the matrix for $\partial_2$ to this nice form.

Consider the chain complex $C_*(\widehat{M})$ associated to the induced 
handle decomposition of the maximal abelian cover $\widehat{M}$ of
$M$.  This is a free 
$\ZH$-chain complex with distinguished basis determined by lifts of 
handles of $M$.  For an appropriate choice of these lifts, we have $\partial_1$ given by $x \mapsto 
x\cdot w$ where $w$ is a column of height $m$ whose $i^{\text{th}}$ entry is 
$h_i -1$.  The map $\partial_2$ is
given by the Alexander-Fox matrix for the presentation $\langle 
x_1,\dots,x_m | r_1 , \dots,r_{m-1} \rangle$ for an appropriate choice 
of the $r_i$'s in their conjugacy classes.  This is an $(m-1 \times m)$ 
matrix whose $(i,j)^{\text{th}}$ entry is given by $\eta (\partial 
r_i/\partial x_j)$ where $\eta$ is the projection $\Z[F] \longrightarrow 
\Z[\pi_1(M,u)] \longrightarrow \ZH$.  Let $e_N$ be an Euler structure
determined by a fundamental family of cells which gives this
``nice'' cellular structure to $\widehat{M}$.  Clearly the $\ZH$--complex for
$\widehat{M}$ must augment to the $\Z$--complex for $M$, hence $\text{aug}(\eta
(\partial r_i/\partial x_j)) = v_{i,j}$.  From before, we know that $v_{i,j}=0$ except in the lower right-hand corner, so $\eta (\partial
r_i/\partial x_j) \in I$ for $i \leq n-1 , j \leq n$.  We claim for $i
\leq n-1 , j \leq n$,
\begin{equation}
|T|\eta (\partial r_i/\partial x_j) =
 -|T|\sum\limits_{p=1}^n \langle k^*_i \cup h^*_j \cup h^*_p ,
 [M] \rangle (h_p - 1) \bmod I^2.
\label{equation-IntCohomFoxAndCup}
\end{equation}
Here $I$ is the augmentation ideal.
We will prove this by looking at $\twid{\eta}$, the composition of $\eta$ with the 
projection $\ZH \longrightarrow \Z [G]$.  If we let $J$ denote 
the augmentation ideal in $\Z [G]$, then it is enough to 
show that for $i \leq n-1 , j \leq n$,
\begin{equation}
\twid{\eta} (\partial r_i/\partial x_j) =
 -\sum\limits_{p=1}^n \langle k^*_i \cup h^*_j \cup h^*_p ,
 [M] \rangle (\twid{h}_p - 1) \bmod J^2.
\label{equation-IntCohomFoxAndTorAug}
\end{equation}
To prove \eqref{equation-IntCohomFoxAndTorAug}, note that $J/J^2$ is isomorphic to the free abelian group $G$ of rank 
$n$ under the map $g\mapsto (g-1) \bmod J^2$, and is thus generated by $\twid{h}_1-1,\dots,\twid{h}_n-1$.  To be precise, for any $g \in
G$, the expansion $g=\prod\limits_{p=1}^n \twid{h}_p^{\langle h_p^*,g\rangle}$ gives
\begin{equation}
\label{equation-I/I^2-basis}
g-1 = \sum\limits_{p=1}^n \langle h_p^* , g \rangle (\twid{h}_p-1)\bmod J^2.
\end{equation}
Also, for any $\alpha \in F, j \leq n$, 
\begin{equation}
\label{equation-aug-Fox-vs-Kron}
\text{aug}(\partial\alpha/\partial
x_j) = \langle h_j^* , \eta(\alpha) \rangle.
\end{equation}
Now $r_i \in [F,F]$
gives an expansion $r_i = \prod\limits_\mu[\alpha_\mu,\beta_\mu]$ a finite
product of commutators in $F$.  Then
\begin{equation*}
\eta(\partial r_i/\partial x_j) = \sum\limits_\mu(\eta(\alpha_\mu) -
1)\eta(\partial\beta_\mu/\partial x_j) +
(1-\eta(\beta_\mu))\eta(\partial\alpha_\mu/\partial x_j).
\end{equation*}
Projecting to $\Z [G]$ we get {\small
\begin{align}
\notag\twid{\eta}(\partial r_i / \partial &x_j)\bmod J^2 =\\
\label{equation-int-cohom-Fox-just-before-surface}&\sum\limits_{p=1}^n \left ( \sum\limits_\mu \langle h_p^* , \eta(\alpha_\mu)
\rangle \langle h_j^* , \eta(\beta_\mu) \rangle - \langle h_p^* ,
\eta(\beta_\mu) \rangle \langle h_j^* , \eta(\alpha_\mu) \rangle
\right ) (\twid{h}_p - 1).
\end{align}}

Now we consider the handlebody $U \subset M$ formed by the (honest) $0$--handle
and the (honest) $1$--handles.  The boundary circle of the core $2$-cell of the $i^{\text{th}}$ honest $2$--handle 
lies in $\partial U$ and represents $r_i$.  The expansion $r_i =
\prod\limits_\mu[\alpha_\mu,\beta_\mu]$ tells us that the circle is the image of the boundary under a map $\sigma_i'\co\Sigma_i'\rightarrow U$, where $\Sigma_i'$ is a surface which have meridians and longitudes with homology classes $a_\mu,b_\mu$ respectively, where $(\sigma_i')_*(a_\mu)=\eta(\alpha_\mu)$ and $(\sigma_i')_*(b_\mu)=\eta(\beta_\mu)$ for all $\mu$.
Let $\Sigma_i$ be $\Sigma_i'$ capped with a disk; the map $\sigma_i'$ extends over the disk to a map $\sigma_i\co\Sigma_i\rightarrow U$ by mapping the capping $2$--disk to the core $2$--cell of the $i^{\text th}$ honest handle.  We may orient $\Sigma_i$ so that the fundamental class $(\sigma_i)_*([\Sigma_i])$ is
represented in $H_2(M)$ by the homology class of the core disk of the
$i^{\text{th}}$ $2$--handle, hence $(\sigma_i)_*([\Sigma_i]) = k_i^* \cap [M]$.  Now
for any $1$-cohomology classes $t_i,t_i'$ of $\Sigma_i$, we have
\begin{equation}
\label{equation-surface-cohomology}
\langle t_i \cup t_i' , [\Sigma_i] \rangle = \sum\limits_\mu \langle t_i ,
a_\mu \rangle \langle t_i' , b_\mu \rangle - \langle t_i ,
b_\mu \rangle \langle t_i' , a_\mu \rangle.
\end{equation}
Pulling $h_j^*$ back to $\Sigma_i$ we get a $1$-cohomology class whose
evaluations on the meridians and longitudes are $\langle h_j^* ,
\eta(\alpha_\mu) \rangle$ and $\langle h_j^* , \eta(\beta_\mu)
\rangle$ respectively.  This proves (everything $\bmod J^2$)
\begin{align}
\label{equation-Int-cohom-Fox-after-surface}
\twid{\eta}(\partial r_i / \partial x_j) & = \sum\limits_{p=1}^n \langle
h^*_p \cup h^*_j , (\sigma_i)_*([\Sigma_i]) \rangle (\twid{h}_p - 1)\\
\notag& = \sum\limits_{p=1}^n \langle h^*_p \cup h^*_j , k^*_i \cap [M] \rangle
(\twid{h}_p - 1)\\
\notag& = \sum\limits_{p=1}^n \langle h^*_p \cup h^*_j \cup k^*_i , [M] \rangle
(\twid{h}_p - 1)\\
\notag& = -\sum\limits_{p=1}^n \langle k^*_i \cup h^*_j \cup h^*_p , [M] \rangle
(\twid{h}_p - 1)\bmod J^2.
\end{align}
This proves \eqref{equation-IntCohomFoxAndTorAug} which in turn proves
\eqref{equation-IntCohomFoxAndCup}. 

Recall by \cite{Tur:3d} II.4.3, we have $\tau(M,e,\omega)\in \ZH$.
We have arranged our handles so that $h_1$ in particular has infinite 
order in $H_1(M)$, so by \eqref{equation:3MbdTorCalc}, we have
\begin{equation*}
(h_1-1)\tau(M,e_N,\omega) = (-1)^{m+1}\tau_0 \det (\Delta(1)).
\end{equation*}
Recall $e_N$ is chosen so that we may use \eqref{equation:3MbdTorCalc}.
We now want to work out $\tau_0$.  For now we work in a very specific
homology basis: 
\[\{[\text{pt}],h_1,\dots,h_n,k^*_1\cap [M],\dots,k^*_{n-1}\cap [M]\}.\]
Let $\twid{\omega}$ denote the homology orientation defined by this basis, and let $[\twid{\omega}/\omega]$ denote the sign of $\twid{\omega}$ with respect to $\omega$.  Later, when we do the $\Det(f)$ calculation, we will use the bases for
$H^1(M)$ and $H^1(M,\partial M)$ given by
$\{h_1^*,\dots,h_n^*\}$ and $\{k_1^*,\dots,k_{n-1}^*\}$
respectively.  Using the basis above (by which we defined $\twid{\omega}$), we compute 
\[
\tau(C_*(M;\R)) = (-1)^{|C_*(M)| + n(m-n)}\det v,
\]
where $v$ is defined as above.  This is a quick
calculation; we may choose our bases of the images of the boundary
maps so that the dimension $2$ and dimension $0$ change of basis matrices
are the identity matrices.  Then the dimension $1$ change of basis
matrix will be the block matrix $\left ( \begin{smallmatrix} 0 &
v \\ {\text{id}}  & 0 \end{smallmatrix} \right )$,
where id represents the $(n\times n)$ identity matrix.  This has
determinant $(-1)^{n(m-n)}\ \det v$.  Another quick calculation
gives $|C_*(M)| = (mn + m + n)\bmod 2$.  Hence $\tau_0 = [\twid{\omega}/\omega]
(-1)^m \sign(\det v)$.  This gives
\begin{equation*}
(h_1 -1)\tau(M,e_N,\omega) = [\twid{\omega}/\omega] (-1)^{m+m+1} \sign(\det v)
\det(\Delta(1)).
\end{equation*}
Let $a$ denote the submatrix of $\Delta$ comprised of the first $n-1$
rows and $n$ columns; thus $a$ is the matrix whose $i,j$ entry is
given by $\eta(\partial r_i/\partial x_j)$ for $1 \leq i \leq n-1$ and
$1\leq j \leq n$.  Let $V$ denote the lower right hand $(m-n \times
m-n)$ matrix $\eta(\partial r_i/\partial x_j)$ for $n \leq i \leq m -
1$ and $n+1 \leq j \leq m$. 
Hence 
\begin{align*}
(h_1-1)\tau(M,e_N,\omega) & = - [\twid{\omega}/\omega]|\det v|\det a(1) \\
& = - [\twid{\omega}/\omega]|T|\ \det a(1)\bmod I^n.
\end{align*}

Define 
\begin{equation*}
\theta_{i,j} = \sum\limits_{p=1}^n \langle k^*_i \cup h^*_j \cup h^*_p ,[M] 
\rangle \twid{h}_p.
\end{equation*}
Then by \eqref{equation-IntCohomFoxAndCup}, the entries of the matrix $|T|a$ modulo $I^2$ can by obtained from the entries of the matrix $|T|\theta$ by replacing each $\twid{h}_p$ by $(h_p-1)$.  From the definition of $q$ above, we see
\begin{equation*}
(h_1-1)\tau(M,e_N,\omega)\bmod I^n = -[\twid{\omega}/\omega]q(\det (\theta(1))).
\end{equation*}

By definition, i.e. \eqref{equation-IntCohomDetDef}, 
\[
\det(\theta(1)) = -[\twid{\omega}/\omega]\Det_{\omega}(f_M)\twid{h}_1.
\]
When we put this together, all of the signs
will neatly cancel out, leaving
\begin{equation*}
(h_1-1)\tau(M,e_N,\omega) \bmod I^n = (h_1-1) q(\Det_{\omega}
(f_M)).
\end{equation*}
Then the map $\bigoplus\limits_{\ell \geq 0}
I^\ell/I^{\ell+1}$ 
defined by $x\in I^\ell/I^{\ell+1}$ maps to $(h_1-1)x\in
I^{\ell+1}/I^{\ell+2}$ is injective, so
\begin{equation*}
\tau(M,e_N,\omega) \bmod I^{n-1} = q(\Det_{\omega} (f_M)).
\end{equation*}
But now recall $\tau(M,e,\omega)$ only differs from
$\tau(M,e_N,\omega)$ by multiplication by an element of $H_1(M)$.  They
are both in $I^{n-2}$, so $\bmod{\ I^{n-1}}$ they are equal.  This
completes the proof.
\end{proof}
\end{subsection}
\end{section}

%
\begin{section}{The Cohomology Ring Mod--$r$}
\label{section-ModCohom}
In this section, we will prove an analogous result to the one in Section~\ref{section-IntCohom}
using cohomology modulo an integer $r$ rather than integral
cohomology.  The integer $r$ will have to be one such that the
first cohomology group with Mod--$r$ coefficients is a free
$\Z_r$--module; for instance if $r$ is prime.  This will also imply
that the first relative cohomology group is a free $\Z_r$--module, so
we will still be able to compute a determinant as in
Section~\ref{section-IntCohom}, however will will need to refine that
determinant slightly.  To do so, we will first introduce the concept
of a {\it paired volume form}, which will play a similar role to the
square volume forms found in \cite{Tur:3d}, III.3

Before anything else, however, let us define the Mod--$r$ torsion.  This is defined when $b_1(M)\geq 2$ so that $\tau(M,e,\omega)\in \Z[H_1(M)]$ for any $e,\omega$.  Then $\tau(M,e,\omega;r)$ is the image of $\tau(M,e,\omega)$ under the projection $\Z[H_1(M)]\rightarrow \Z_r[H_1(M)]$ induced by the coefficient projection $\Z\rightarrow \Z_r$.  Note that if $r=p_1^{e_1}\cdot p_2^{e_2}\cdots p_k^{e_k}$ where $p_1,\dots ,p_k$ are primes, then $\Z_r[H_1(M)]$ splits naturally as $\Z_{p_1^{e_1}}[H_1(M)]\oplus\dots\oplus\Z_{p_k^{e_k}}[H_1(M)]$ and $\tau(M,e,\omega;r)$ splits as $\tau(M,e,\omega;p_1^{e_1}) + \dots + \tau(M,e,\omega;p_k^{e_k})$, so understanding Mod--$r$ torsion when $r$ is a power of a prime is sufficient to understand it for any $r$.

One may also define the Mod--$r$ torsion when $b_1(M)=1$ by using Turaev's ``polynomial part'' $[\tau]$ of the torsion; see \cite{Tur:3d}, II.3.  Theorem~\ref{theorem-mod_r-vs-cohom} is true in this case as well, and one can use the argument in \cite{Tur:3d} Theorem~III.4.3 when the first Betti number is $1$ (the last paragraph of the proof).  
%
\begin{subsection}{Determinants}
\label{subsection-ModCohom-Dets}
\begin{subsubsection}{Volume Forms}
\label{subsubsection-MixedVolumeForms}
First we recall some definitions from \cite{Tur:3d}, III.3.  If $N$ is
a finite rank free module over $R$, a commutative ring with $1$, then a
{\it volume form} $\omega$ on $N$ is a map which assigns to each basis
$a$ of $N$ a scalar $\omega(a)\in R$ such that $\omega(a) =
[a/b]\omega(b)$ for any bases $a,b$.  A {\it square volume form} is a
map $\Omega$ which also assigns a scalar to each basis, but the change
of basis formula is $\Omega(a)=[a/b]^2\Omega(b)$.  Naturally, the
square of a volume form is a square volume form.  This notion is useful when working with closed manifolds as in \cite{Tur:3d}, III.3, but we must use a slightly different form in the case of a nonempty boundary, though in the same spirit.  
If $K,L$ are two
finite rank free $R$--modules, then a {\it paired volume form} on
$K\times L$ is a map $\mu$ from (ordered) pairs of bases of $K$
and $L$ to $R$ such that $\mu(a',b')=[a'/a][b'/b]\mu(a,b)$ where
$a,a'$ are bases of $K$ and $b,b'$ are bases of $L$. Note that the
product of a volume form on $K$ with a volume form on $L$ is a paired
volume form on $K\times L$, so this notion is very similar to the
notion of a square volume form.  We say a paired volume form is {\it
non-degenerate} if its image lies in the units of $R$, or equivalently
if there is a basis $a$ of $K$ and a basis $b$ of $L$ so that
$\mu(a,b)=1$.  Note that we may easily construct a non-degenerate
paired volume form given distinguished bases $a,b$ of $K,L$
respectively by {\it assigning} $\mu(a,b)=1$, and extending to other
bases by the change of basis formula.

Note the following properties of paired volume forms:
\begin{enumerate}
\item\label{enumitem-volForms-bilin} If $B\co K\times L\rightarrow R$ is a bilinear form, where $K$ and $L$ are isomorphic $R$--modules, then $\mu(a,b)=\det(B_{a,b})$ is a paired volume form, where $B_{a,b}$ is the matrix of $B$ with respect to the bases $a$ and $b$ of $K$ and $L$ respectively, and $\mu$ is non-degenerate if and only if $B$ is a nondegenerate form, i.e. if $B$ induces an isomorphism $K\rightarrow \Hom_R(L,R)$.
\item\label{enumitem-volForms-freeZ} If $K,L$ are free $\Z$--modules of finite rank $r_K$ and $r_L$ respectively, and $\omega$ is an orientation on $(K\times L)\otimes \R$, then there is a non-degenerate paired volume form $\mu_\omega$ on $K\times L$ such that $\mu_\omega(a,b) = 1$ if the basis $a_1\otimes 1,\dots,a_{r_K}\otimes 1,b_1\otimes 1,\dots,b_{r_L}\otimes 1$ is positively oriented with respect to $\omega$ (and obviously $\mu_\omega$ assigns $-1$ to bases which are negatively oriented with respect to $\omega$).
\item\label{enumitem-volForms-SES} If $0\longrightarrow K_1\longrightarrow K\longrightarrow K_2\longrightarrow 0$ and $0\longrightarrow L_1\longrightarrow L\longrightarrow L_2\longrightarrow 0$ are short exact sequences of finite rank free $R$--modules and $\mu_1,\mu_2$ are paired volume forms on $K_1\times L_1,K_2\times L_2$, then there is an induced paired volume form on $K\times L$, which is non-degenerate if and only if $\mu_1$ and $\mu_2$ are both non-degenerate.  To construct this, let $a_i,b_i$ be bases of $K_i,L_i$ respectively.  Then we can construct the bases $a_1a_2$ and $b_1b_2$ of $K$ and $L$ respectively by concatenating the image of the basis $a_1$ with a lift of the basis $a_2$ in $K$, and similarly for $b_1b_2$ in $L$.  Then for any bases $a$ and $b$ of $K$ and $L$, define $\mu(a,b)=[a/a_1a_2][b/b_1b_2]\mu_1(a_1,b_1)\mu_2(a_2,b_2)$.  
\item\label{enumitem-volForms-dual} A non-degenerate paired volume form $\mu$ on $K\times L$ induces a non-degenerate paired volume form $\mu^*$ on $K^*\times L^*\iso (K\times L)^* =\Hom_R(K\times L , R)$ by $\mu^*(a^*,b^*) = (\mu(a,b))^{-1}$ where $a^*$ is the basis of $K^*$ dual to a basis $a$ of $K$, and similarly for $b,b^*$.  
\item\label{enumitem-volForms-quotient} If $\phi\co R\rightarrow S$ is a surjection of rings, and $\mu$ is a nondegenerate paired volume form on the free $R$--modules $K\times L$, then there is an induced paired volume form $\mu_\phi$ on $K\otimes_S S \times L\otimes_S S$ given by $\mu_\phi(a\otimes 1,b\otimes 1)=1$ if $a,b$ are bases of $K,L$ such that $\mu(a,b)=1$. 
\end{enumerate}
\end{subsubsection}
\begin{subsubsection}{The Refined Determinant}
\label{subsubsection-ModCohomRefinedDet}
Now given free $R$--modules $K,L$ of finite ranks $n$ and $n-1$ respectively ($n\geq 2$), and given $f\co L\times K\times K\rightarrow R$ an $R$--map as in Lemma~\ref{lemma-IntCohomDet}, and given a paired volume form $\mu$ on $K\times L$, we can construct the $\mu$--refined determinant, $\Det_\mu(f)$, to be 
\begin{equation}
\label{equation-volFormRefinedDet}
\Det_\mu(f)=\mu^*(a^*,b^*)d(f,a,b)
\end{equation} 
where $d$ is defined as in Lemma~\ref{lemma-IntCohomDet}.  We can define this for any bases $a,b$ of $K,L$ respectively (and $a^*$, $b^*$ the dual bases as usual), but by the properties of $d$ and $\mu$, this is independent of the chosen bases.  Note that this will simply be the determinant taken with respect to any bases $a,b$ with $\mu(a,b)=1$ if such bases exist.
\end{subsubsection}
\begin{subsubsection}{Constructing Paired Volume Forms}
\label{subsubsection-ModCohomConstructForm}
We now construct a paired volume form which we will later use to refine the determinant for mod--$r$ cohomology.  We will construct this form in pieces, and assemble them via methods enumerated above.  Let $H,H'$ be finite abelian groups which are isomorphic, though we will not fix a particular isomorphism.  (These groups will appear later as the torsion groups $\Tors(H_1(M))$ and $\Tors(H_1(M,\partial M))$).  Let $p\geq 2$ be a prime integer dividing $\vert H \vert$.  Let $r=p^s$ for some $s\geq 1$ such that $H/r$ is a direct sum of copies of $\Z_r$, so that we can think of $H/r$ as a finite rank free $\Z_r$--module (and similarly for $H'/r$, since $H,H'$ are isomorphic).  We will now show how to construct a paired volume form on $H/r\times H'/r$ from a bilinear form $L\co H\times H'\rightarrow\Q/\Z$.  First, we repeat some definitions from \cite{Tur:3d}.

Let $H_{(p)}$ be the subgroup of $H$ consisting of all elements annihilated by a power of $p$ (similarly for $H'_{(p)}$).  A sequence $h=(h_1,\dots ,h_n)$ of nonzero elements of $H_{(p)}$ is a {\it pseudo-basis} if $H_{(p)}$ is a direct sum of the cyclic subgroups generated by $h_1,\dots ,h_n$ and the order of $h_i$ in $H$ is less than or equal to the order of $h_j$ for $i\leq j$.  In other words, if the order of $h_i$ is $p^{s_i}$, with $s_i\geq 1$, then $s_1\leq s_2\leq\dots\leq s_n$.  This sequence $(s_1,\dots ,s_n)$ is determined by $H_{(p)}$ and does not depend on the choice of pseudo-basis $h$.  Note $s\leq s_1$ since if we have a summand of order $p^k$ for $k < s$, then projecting to $H/r$ there is still a summand of order $p^k$, which contradicts our assumption that $H/r$ is a sum of several $\Z_r$'s.  Projecting a pseudo-basis to $H_{(p)}/r=H/r$ we get a basis $\overline{h}$ of the $\Z_r$--module $H/r$.

Let $L\co H\times H'\rightarrow \Q/\Z$ be a bilinear form.  We will say $L$ is nondegenerate if the map induced by $L$ from $H\rightarrow \Hom_\Z(H',\Q/\Z)$ is an isomorphism (since all of the groups involved are finite and of the same order, this is equivalent to the map being an injection or a surjection).  Note if $z'\in H'_{(p)}$, then $z'$ has order $p^k$ for some $k\geq s$.  Then for any $z\in H , L(z,z')\in (p^{-k}\Z)/\Z$, and therefore $p^kL(z,z')\in \Z/(p^k\Z)$.  Projecting this to $\Z_r$, we obtain an element which we will call $z\cdot z'$.  Note we can do something similar if $z$ has order $p^k$ and $z'$ does not necessarily, and that they clearly agree if both $z,z'$ have order a power of $p$.  We have now constructed $z\cdot z'$, a $\Z_r$ pairing on $H\times H'$.  The following is the analogue of Lemma~III.3.4.1 in \cite{Tur:3d} (the proof is a direct generalization of the proof found there as well).

\begin{lemma}
\label{lemma-ModCohomConstructPair}
There is a unique paired volume form $\mu_L^r$ on $H/r\times H'/r$ such that for any pseudo-bases $h=(h_1,\dots ,h_n), k=(k_1,\dots ,k_n)$ of $H_{(p)}, H'_{(p)}$ respectively,
\begin{equation}
\label{equation-ModCohomConstructPair}
\mu_L^r(\overline{h},\overline{k}) = \det(h_i\cdot k_j) \in \Z_r.
\end{equation}
Also, if $L$ is nondegenerate, then so is $\mu_L^r$.
\end{lemma}
\begin{proof}
It is clear that given distinguished pseudo-bases $h,k$ then we can construct a paired volume form $\mu_{(h,k)}$ by $\mu_{(h,k)}(a,b)=[a/\overline{h}][b/\overline{k}] \det(h_i\cdot k_j)$ for any bases $a,b$ of $H/r, H'/r$ respectively. Then $\mu_{(h,k)}(\overline{h},\overline{k})=\det(h_i\cdot k_j)$, so we would like to define $\mu_L^r = \mu_{(h,k)}$.  We now prove that the definition of $\mu_{(h,k)}$ does not actually depend on $h$ or $k$.  To prove this, it suffices to show that for any other pseudo-bases $x=(x_1,\dots , x_n)$ of $H_{(p)}$ and $y=(y_1,\dots ,y_n)$ of $H'_{(p)}$,
\begin{equation}
\label{equation-ModCohomConstructPairIndepend}
\det(x_i\cdot y_j)=[\overline{x}/\overline{h}][\overline{y}/\overline{k}]\det(h_i\cdot k_j).
\end{equation}
By symmetry, to prove \eqref{equation-ModCohomConstructPairIndepend} we need only show
\begin{equation}
\label{equation-ModCohomConstructPairIndependSingle}
\det(x_i\cdot k_j) = [\overline{x}/\overline{h}]\det(h_i\cdot k_j).
\end{equation}
Now $x$ is a pseudo-basis for $H_{(p)}$, so the order of $x_i$ is equal to the order of $h_i$ for each $i$.  It is clear that if $x$ is just a permutation of $h$ (the permutation can only permute elements of the same order), then the basis $\overline{x}$ of $H/r$ is the same permutation of the basis $\overline{h}$, and then \eqref{equation-ModCohomConstructPairIndependSingle} is clear.  So now, we may assume that each $x_i$ generates the same cyclic subgroup of $H_{(p)}$ as the corresponding $h_i$.  Then for each $i$, there is some $c_i\in\Z$, with $c_i$ coprime to $p^{s_i}$, hence coprime to $r=p^s$ (in fact, coprime to $p$), with $x_i=c_i h_i$.  But then $x_i\cdot k_j = (c_i \pmod r ) h_i\cdot k_j$, so $\det(x_i\cdot k_j) = \prod_i (c_i\pmod r)\det(h_i\cdot k_j)$.  But clearly $[\overline{x}/\overline{h}] = \prod_i c_i\pmod r$, so the proof of \eqref{equation-ModCohomConstructPairIndependSingle} is completed, and \eqref{equation-ModCohomConstructPairIndepend} follows.

Now if $L$ is nondegenerate, then to show that $\mu_L^r$ is nondegenerate, we must simply show that $\det(h_i\cdot k_j)\in \Z_r^\times$ for any pseudo-bases $h,k$ of $H_{(p)},H'_{(p)}$ respectively.  
$L$ nondegenerate means that the map induced by $L$, $\twid{L}\co H\rightarrow \Hom_\Z(H',\Q/\Z),$ is a bijection.  Then, in particular, the restriction of $\twid{L}$ to $H_{(p)}$ is also bijective on its image $\Hom_\Z(H'_{(p)} , \Q/\Z)$.   This means, for $k$ any pseudo-basis of $H'_{(p)}$, for each $k_j$ there is an $x_j\in H_{(p)}$ with $L(x_i,k_j)=\delta_{i,j}p^{-s_j}$, i.e. $x_i\cdot k_j =\delta_{i,j}$.  Then \eqref{equation-ModCohomConstructPairIndependSingle} gives us our result, that $\det(h_i\cdot k_j)\in \Z_r^\times$ for any pseudo-bases $h,k$ of $H_{(p)},H'_{(p)}$ respectively.
\end{proof}
\end{subsubsection}
\begin{subsubsection}{The $\Q/\Z$ linking form}
\label{subsubsection-ModCohom_Q/Z_linking}
We now construct a linking form on 
\[
\Tors(H_1(M))\times \Tors(H_1(M,\partial M)).
\]
We use a slightly different construction from the one in \cite{Tur:3d}, though one may easily verify that the end results are the same.

From the Universal Coefficient Theorem, there is an exact sequence
\[ 0 \rightarrow H_2(M)\otimes \Q/\Z \rightarrow H_2(M;\Q/\Z) \rightarrow \Tor(H_1(M),\Q/\Z) \rightarrow 0\]
but $\Tors(H_1(M))$ is canonically isomorphic to $\Tor(H_1(M),\Q/\Z)$ by 
\begin{align*}
\Tors(H_1(M))&\iso\Tors(H_1(M))\otimes\Z\\
&\iso\Tor(\Tors(H_1(M)) , \Q/\Z)\\
&\iso\Tor(H_1(M) , \Q/\Z)
\end{align*}
With this in mind, our exact sequence becomes
\[ 0 \rightarrow H_2(M)\otimes \Q/\Z \rightarrow H_2(M;\Q/\Z) \rightarrow \Tors(H_1(M)) \rightarrow 0.\]
Now choose elements $a\in\Tors(H_1(M))$ and $b\in\Tors(H_1(M,\partial M))$; we want to define their linking $L_M(a,b)\in\Q/\Z$.  So choose $\overline{a}\in H_2(M;\Q/\Z)$ mapping to $a$, then let $\alpha\in H^1(M,\partial M;\Q/\Z)$ be Poincar\'{e} dual to $\overline{a}$, i.e. $\alpha \cap [M]=\overline{a}$.  Then we define $L_M(a,b)=\langle \alpha , b \rangle \in \Q/\Z$.  One may show that starting with the sequence for $H_2(M,\partial M;\Q/\Z)$ instead of the sequence for $H_2(M;\Q/\Z)$ yields the same form.
\end{subsubsection}
\begin{subsubsection}{Constructing the Paired Volume Form for Cohomology}
\label{subsubsection-ConstructingPairedFormModCohom}
Let $\omega$ be a homology orientation.  Then we have split exact sequences{\small
\[ 0 \rightarrow \Tors(H_1(M))\rightarrow H_1(M)\rightarrow H_1(M)/\Tors(H_1(M))\rightarrow 0,\]}{\footnotesize
\[ 0 \rightarrow \Tors(H_1(M,\partial M))\rightarrow H_1(M,\partial M)\rightarrow H_1(M,\partial M)/\Tors(H_1(M,\partial M))\rightarrow 0.\]}
Both of these sequences split, so they also split modulo $r$, and note
\[
H_1(M)/r\iso H_1(M;\Z_r) \text{\ and\ } H_1(M,\partial M)/r\iso H_1(M,\partial M;\Z_r).
\]
The homology orientation induces a nondegenerate paired volume form on the free $\Z$--module \[H_1(M)/\Tors(H_1(M))\times H_1(M,\partial M)/\Tors(H_1(M,\partial M))\] which induces a nondegenerate paired volume form on \[(H_1(M)/\Tors(H_1(M)))/r\times (H_1(M,\partial M)/\Tors(H_1(M,\partial M)))/r.\]  Above, we showed how to construct a nondegenerate paired volume form (induced by the $\Q/\Z$--linking form) on \[\Tors(H_1(M))/r\times \Tors(H_1(M,\partial M))/r.\]  We may combine these to give a nondegenerate paired volume form on $H_1(M;\Z_r)\times H_1(M,\partial M;\Z_r)$, which in turn gives us a nondegenerate paired volume form on the duals.  We will denote the canonical Mod--$r$ paired volume form by on $H_1(M;\Z_r)\times H_1(M,\partial M;\Z_r)$ by $\mu^r_M$ and the refined determinant of the form $f^r_M$ on $H^1(M,\partial M;\Z_r)\times H^1(M;\Z_r)\times H^1(M;\Z_r)$ by $\Det_r(f^r_M)$.
\end{subsubsection}
\end{subsection}
%
\begin{subsection}{Relationship to Torsion}
\label{subsection-ModCohom-Relns}
We will now let $I$ denote the augmentation ideal of $\Z_r[H_1(M)]$ instead of the augmentation ideal of $\Z[H_1(M)]$ (the augmentation ideal of $\Z_r[H_1(M)]$ is the image of the augmentation ideal of $\Z[H_1(M)]$ under the map induced by the coefficient projection $\Z\rightarrow \Z_r$).  We now recall a definition from \cite{Tur:3d} - we define $q_r\co S(H_1(M)/r)\rightarrow \bigoplus\limits_{\ell\geq 0}I^\ell/I^{\ell + 1}$ by
\begin{equation}
\label{equation-q_r-def}
q_r(g_1,\dots ,g_\ell) = \prod\limits_{i=1}^\ell (\twid{g_i}-1) \pmod{I^{\ell +1}}
\end{equation}
where $\twid{g_i}$ is a lift of $g_i$ to $H_1(M)$ (the proof that this is independent of the lift is in \cite{Tur:3d}).

Before we state the main theorem, we briefly discuss Mod--$r$ surfaces.  In particular, we give equivalent equations to \eqref{equation-surface-cohomology}.  An equivalent definition of Mod--$r$ surfaces can be found in \cite{Tur:3d} Section~XII.3.

\begin{subsubsection}{Mod--$r$ surfaces}
\label{subsubsection-Mod-r-surfaces}
Let $G(\mathcal{M},\mathcal{N};r)$ be a group with generators $\alpha_\mu,\beta_\mu,\gamma_\nu$ where $\mu,\nu$ run over finite indexing sets $\mathcal{M},\mathcal{N}$ respectively, with a single relator $\rho=\prod\limits_\mu [\alpha_\mu,\beta_\mu]\prod\limits_\nu \gamma_\nu^r$.  Let $X(\mathcal{M},\mathcal{N};r)$ be a connected CW--complex with a single $0$--cell, $1$--cells $a_\mu,b_\mu,c_\nu$ (so that we can consider $\pi_1(X(\mathcal{M},\mathcal{N};r))$ to be generated by $\alpha_\mu,\beta_\mu,\gamma_\nu$), and a single $2$--cell attached along $\rho$, so it is obvious that \[\pi_1(X(\mathcal{M},\mathcal{N};r))\iso G(\mathcal{M},\mathcal{N};r).\]  Then $H_2(X(\mathcal{M},\mathcal{N};r);\Z_r)\iso \Z_r$, so let $[X(\mathcal{M},\mathcal{N};r)]$ be the generator of $H_2(X(\mathcal{M},\mathcal{N};r);\Z_r)$ given by the homology class of the two cell (whose boundary is zero modulo $r$).  Now if $t,t'\in H^1(X(\mathcal{M},\mathcal{N};r);\Z_r)$, then let us compute \[\langle t\cup t' ,[X(\mathcal{M},\mathcal{N};r)]\rangle = \varepsilon_*\left( (t\cup t')\cap [X(\mathcal{M},\mathcal{N};r)]\right)\] where $\varepsilon_*\co H_0(X(\mathcal{M},\mathcal{N};r);\Z_r)\rightarrow \Z_r$ is simply augmentation, $[\text{pt}]\mapsto 1$.  Let $a_\mu,b_\mu,c_\nu\in H_1(X(\mathcal{M},\mathcal{N};r);\Z_r)$ be the homology classes modulo $r$ of $\alpha_\mu,\beta_\mu,\gamma_\nu$ respectively.  We now claim
\begin{lemma}
\label{lemma-mod-r-surface-cohomology}
Let $t,t'\in H^1(X(\mathcal{M},\mathcal{N};r);\Z_r)$.  If $r$ is odd, then 
\begin{equation}
\label{equation-odd-r-surface-cohomology}
\langle t\cup t' , [X(\mathcal{M},\mathcal{N};r)]\rangle = \sum\limits_\mu \langle t , a_\mu\rangle\langle t',b_\mu\rangle - \langle t , b_\mu\rangle\langle t',a_\mu\rangle .
\end{equation}
If $r$ is even, then
\begin{equation}
\label{equation-even-r-surface-cohomology}
\text{\footnotesize $\langle t\cup t' , [X(\mathcal{M},\mathcal{N};r)]\rangle = \sum\limits_\mu \langle t , a_\mu\rangle\langle t',b_\mu\rangle - \langle t , b_\mu\rangle\langle t',a_\mu\rangle + \frac{r}{2}\;\sum\limits_\nu \langle t , c_\nu\rangle\langle t' , c_\nu\rangle$} .
\end{equation}
\end{lemma}
\begin{proof}
If we let $a_\mu^*,b_\mu^*,c_\nu^*\in H^1(X(\mathcal{M},\mathcal{N};r);\Z_r)$ be dual to $a_\mu,b_\mu,c_\nu$ under $\langle\cdot , \cdot\rangle$, then $1 = \langle a^*_\mu\cup b^*_\mu , [X(\mathcal{M},\mathcal{N};r)]\rangle = -\langle b^*_\mu\cup a^*_\mu , [X(\mathcal{M},\mathcal{N};r)]\rangle$.  Clearly $c_\nu\cup c_\nu$ is $2$--torsion for any $r$, and one can also show that all other cup products are zero (this follows from induction and a relatively simple Mayer--Vietoris argument).  So the claim for $r$ odd is completed.  By the same Mayer--Vietoris argument, for even $r$, we only need to show the statement for $\mathcal{M}$ empty, and $\mathcal{N}$ only having one element, i.e. for even $r$, and a CW complex $X$ with one $0$--cell, one $1$--cell $c$, and one $2$--cell with boundary $r\cdot c$, we need to show $\langle c^2,[X]\rangle = \frac{r}{2}$.  But this follows from simply noting that $X$ is the $2$--skeleton of a $K(\Z_r,1)$.  A more complete proof may be found in \cite{Hatcher:AT} Chapter~3, Example~3.9.
\end{proof}
\end{subsubsection}

We are now ready to state the main theorem of this section.
\begin{theorem}
\label{theorem-mod_r-vs-cohom}
Let $r$ be a power of a prime such that $H_1(M)/r=H_1(M;\Z_r)$ is a free $\Z_r$--module of rank $b\geq 2$.  Let \[T=|\Tors(H_1(M))|/r.\]  Then for any Euler structure $e$ and homology orientation $\omega$, \[\tau(M,e,\omega;r)\in I^{b-2},\] and
\begin{equation}
\label{equation-mod_r-vs-cohom}
\tau(M,e,\omega;r) = T\cdot q_r(\Det_r(f^r_M)) \pmod{I^{b-1}}.
\end{equation}
\end{theorem}
As in Theorem~\ref{theorem-IntCohom-vs-tors}, that $\tau(M,e,\omega;r)\in I^{b-2}$ is proved in \cite{Tur:3d} II.4.4, we are more concerned with the residue class modulo $I^{b-1}$.
\begin{proof}
The proof is similar to the proof of Theorem~\ref{theorem-IntCohom-vs-tors}, and again is the method of the proof of \cite{Tur:3d} Theorem~4.3 with modifications to apply it to manifolds with nonvoid boundary.  Suppose $r=p^s$, where $p\geq 2$ is prime and $s\geq 1$.  Let $n=b_1(M)$.  Then $H_1(M)$ splits as $\Z^n\times (\Tors(H_1(M)))_{(p)}\times H'$ and $H_1(M,\partial M)$ splits as $\Z^{n-1}\times (\Tors(H_1(M,\partial M)))_{(p)}\times H''$ where the subscript of $(p)$ denotes the maximal subgroup of a finite group whose order is a power of $p$ and $H',H''$ are (isomorphic) subgroups of $\Tors(H_1(M))$ and $\Tors(H_1(M,\partial M))$ respectively with $|H'|=|H''|=T$.  We again choose a handle decomposition of $M$ and the dual relative handle decomposition of $(M,\partial M)$ with $1$ honest $0$--handle, $m$ honest $1$--handles, $m-1$ honest $2$--handles, and no other handles, where $m\geq b\geq n$.  Let $x_1,\dots ,x_m\in \pi_1(M)$ be the generators of $\pi_1(M)$ (based at the $0$--cell) given by the core $1$--cells of the honest $1$--handles, and let $h_1,\dots h_m$ denote their homology classes.  Let $k_1,\dots,k_{m-1}$ denote the classes in $H_1(M,\partial M)$ of the core cells of the relative $1$--handles, and let $r_1,\dots ,r_{m-1}$ be the relators in $F=\langle x_1,\dots ,x_m\rangle$ given by the attaching maps of the honest $2$--cells.  Now, as in the proof of Theorem~\ref{theorem-IntCohom-vs-tors}, we want to rearrange handles for a more convenient decomposition.

As in the proof of Theorem~\ref{theorem-IntCohom-vs-tors}, we can arrange the handle decomposition (by sliding handles) so that $h_1,\dots,h_n$ are generators modulo $\Tors(H_1(M))$ and $h_{n+1}\dots,h_m\in \Tors(H_1(M))$.  We can also arrange for $h_{n+1},\dots,h_b$ to be a pseudo-basis of $(\Tors(H_1(M)))_{(p)}$ by essentially the same method of sliding handles.  The last handles then have the property $h_{b+1},\dots,h_m\in H'$.  Let $p^{s_1},\dots,p^{s_{b-n}}$ be the orders of $h_{n+1},\dots,h_b$ respectively, and we may assume $s_1\leq s_2\leq\dots\leq s_{b-n}\leq s$.

Now we will denote by $\twid{h}$ the projection of $h\in H_1(M)$ to $H_1(M)/r$, then $\twid{h}_1,\dots,\twid{h}_b$ is a basis for $H_1(M)/r$ over $\Z_r$ and $\twid{h}_i=1$ for $i>b$.  Let $h_i^*\in H^1(M;\Z_r)$ for $i \leq b$ such that $\langle h_i^*,\twid{h}_j\rangle =\delta_{i,j}$.

Let $k_i$ denote the class in $H_1(M,\partial M)$ of the $i^\text{th}$ relative handle, using the methods in the proof of Theorem~\ref{theorem-IntCohom-vs-tors} and the methods above, we can arrange so that $k_1,\dots,k_{n-1}$ are generators modulo $\Tors(H_1(M,\partial M))$, $k_{n},\dots,k_{b-1}$ form a pseudo-basis of $(\Tors(H_1(M,\partial M)))_{(p)}$ (they also have orders $p^{s_1}, \dots, p^{s_{b-n}}$) and $k_b,\dots,k_{m-1}\in T''$.  This means, using $\twid{k}$ to denote projection of $k\in H_1(M,\partial M)$ to $H_1(M,\partial M)/r$, that $\twid{k}_1,\dots,\twid{k}_{b-1}$ is a basis for $H_1(M,\partial M)/r$ over $\Z_r$ and $\twid{k}_i=1$ for $i>b-1$.  As above, let $k_i^*\in H^1(M,\partial M;\Z_r)$ for $i\leq b-1$ be such that $\langle k_i , \twid{k}_j \rangle = \delta_{i,j}$.

The matrix for the boundary map from dimension two to dimension one in $C_*(M)$ decomposes, as in the proof of Theorem~\ref{theorem-IntCohom-vs-tors}, as $\left(\begin{smallmatrix}0&0\\0&v\end{smallmatrix}\right)$, where $v$ is a square presentation matrix for $\Tors(H_1(M))$.  Using the handle decomposition above, $v$ can be split as the direct sum of a diagonal matrix (with $p^{s_1},\dots,p^{s_{b-n}}$ along the diagonal) and a square matrix $v'$ which is a presentation matrix for $H'$ (and its transpose a presentation matrix for $H''$), hence $\det(v')=\pm T$.  The diagonal submatrix of $v$ (consisting of powers of $p$) arises as follows:  for $n+1\leq i\leq b$, $h_i$ has order $p^{s_{i-n}}$ (according to the above argument); in fact, $p^{s_{i-n}}h_i$ is the boundary of the $2$--cell transverse to $k_i$.  This $2$--cell has boundary zero in $\Q/\Z$, and its homology class in $H_2(M;\Q/\Z)$ is Poincar\'{e} dual to the class of $k_i^*$ in $H^1(M,\partial M;\Q/\Z)$, which is dual under evaluation to the class of $k_i$ in $\Tors(H_1(M,\partial M))$.  This process is the precise process used in the construction of the linking pairing, first lifting an element of $\Tors(H_1(M))$ to $H_2(M;\Q/\Z)$ and then using Poincar\'{e} duality to get an element dual (under evaluation) to an element of $\Tors(H_1(M,\partial M))$.  The fact that this matrix is diagonal with $p^{s_{i-n}}$ running down the diagonal means that for $n+1\leq i\leq b$ and $n\leq j\leq b-1$,
\begin{equation}
\label{equation-h_i,k_j-linking}
(h_i\cdot k_j)=\delta_{i,j}.
\end{equation}

Now as in the proof of Theorem~\ref{theorem-IntCohom-vs-tors}, $r_1,\dots,r_{n-1}\in [F,F]$, and the above argument shows $r_n,\dots,r_{b-1}$ can each be expanded as $r_i=\prod\limits_{\mu\in\mathcal{M}_i} [\alpha_\mu,\beta_\mu]\prod\limits_{\nu\in\mathcal{N}_i} \gamma^r_\nu$, so we need to use Lemma~\ref{lemma-mod-r-surface-cohomology}.  Henceforth, we will suppress the $\mathcal{M}_i,\mathcal{N}_i$ notation for simplicity.

Now let $\text{pr}\co\Z[H_1(M)]\rightarrow \Z_r[H_1(M)]$ be coefficient projection, let $\eta\co\Z[F]\rightarrow\Z[H_1(M)]$ be induced by the projection $F\rightarrow H_1(M)$ (through $\pi_1(M)$), and let $p\co\Z_r[H_1(M)]\rightarrow \Z_r[H_1(M)/r]$ be induced by $H_1(M)\rightarrow H_1(M)/r$.  Finally, we will also denote by $\eta_r=p\circ\text{pr}\circ\eta$.  We now prove the analogue of \eqref{equation-IntCohomFoxAndCup} which is, for $i\leq b-1, j\leq b$
\begin{equation}
\label{equation-Mod-r-CohomFoxAndCup}
(\text{pr}\circ\eta)(\partial r_i/\partial x_j)=-\sum\limits_{p=1}^b f^r_M(k_i^*,h_j^*,h_p^*)(h_p-1) \pmod{I^2}.
\end{equation}
We will prove \eqref{equation-Mod-r-CohomFoxAndCup} by proving the analogue of \eqref{equation-IntCohomFoxAndTorAug}, which is
\begin{equation}
\label{equation-Mod-r-CohomFoxAndCupAug}
\eta_r(\partial r_i/\partial x_j)=-\sum\limits_{p=1}^b f^r_M(k_i^*,h_j^*,h_p^*)(\widetilde{h}_p-1) \pmod{J^2}.
\end{equation}
Note \eqref{equation-Mod-r-CohomFoxAndCup} follows from \eqref{equation-Mod-r-CohomFoxAndCupAug} since $p$ induces an isomorphism \[\Z_r[H_1(M)]/I^2\rightarrow \Z_r[H_1(M)/r]/J^2\] (where $J$ is the augmentation ideal in $\Z_r[H_1(M)/r]$).  This follows from noting for any $h\in H_1(M)$, $h^r-1\in I^2$ since $(h^r-1)=(h-1)(1+h+\dots +h^{r-1})$, and $(1+h+\dots +h^{r-1})=(h-1)+(h^2-1)+\dots +(h^{r-1}-1)$ in $\Z_r[H_1(M)]$.  To prove \eqref{equation-Mod-r-CohomFoxAndCupAug}, we need to note that \eqref{equation-I/I^2-basis} and \eqref{equation-aug-Fox-vs-Kron} can be used here {\it mutatis mutandis}; indeed, for $c\in H_1(M)/r$, we may use the same formula as \eqref{equation-I/I^2-basis} with slightly different meaning to the symbols
\begin{equation}
\label{equation-Mod-r-J/J^2-basis}
c-1=\sum\limits_{p=1}^b \langle h_p^*,c\rangle (\twid{h}_p -1) \pmod{J^2}.
\end{equation}
Also, for any $\alpha\in F,j\leq b$, if we let $\text{aug}_r$ denote $\text{aug}\circ p\circ\text{pr}\circ\eta$, $\text{aug}_r\co\Z[F]\rightarrow \Z_r$,
\begin{equation}
\label{equation-Mod-r-aug-Fox-vs-Kron}
\text{aug}_r(\partial\alpha/\partial x_j)=\langle h_j^* , \eta_r(\alpha)\rangle.
\end{equation}

For $1\leq i\leq n-1$, we may compute $\partial r_i/\partial x_j$ by
\begin{equation*}
\eta(\partial r_i/\partial x_j) = \sum\limits_\mu(\eta(\alpha_\mu) -
1)\eta(\partial\beta_\mu/\partial x_j) +
(1-\eta(\beta_\mu))\eta(\partial\alpha_\mu/\partial x_j).
\end{equation*}
For $n\leq i\leq b-1$, we must add a term for the $\gamma_\nu$'s
\begin{align}
\notag\eta(\partial r_i/\partial x_j) = &\sum\limits_\mu\big( (\eta(\alpha_\mu) -
1)\eta(\partial\beta_\mu/\partial x_j) +
(1-\eta(\beta_\mu))\eta(\partial\alpha_\mu/\partial x_j)\big)\\
\label{equation-Mod-r-Fox-extra-gamma}&+ \sum\limits_\nu\left( \eta(\partial\gamma_\nu/\partial x_j)(1+\gamma_\nu+\dots+\gamma_\nu^{r-1})\right) .
\end{align}
For any $r$, any $c\in H_1(M)$, working modulo $I^2$, 
\begin{align*}
\sum\limits_{\ell=0}^{r-1} c^\ell &= \sum\limits_{\ell=0}^{r-1} \left( 1+(c-1)\right)^\ell\\
&= \sum\limits_{\ell=0}^{r-1}\sum\limits_{s=0}^\ell \binom{\ell}{s} 1^{\ell-s} (c-1)^{s}\\
&= \sum\limits_{\ell=0}^{r-1} 1 + \ell (c-1) \pmod{I^2}\\
&= \sum\limits_{\ell=0}^{r-1} \ell (c-1)\\
&= (c-1)r(r-1)/2.
\end{align*}
So if $r$ is odd, then each extra $\gamma_\nu$ term is zero modulo $I^2$ in \eqref{equation-Mod-r-Fox-extra-gamma}, and if $r$ is even (in our case, a power of two), then for each $\nu$ (applying \eqref{equation-Mod-r-J/J^2-basis}), 
\begin{align}
\notag1+\eta_r(\gamma_\nu)+\dots+\eta_r(\gamma_\nu)^{r-1}&=-\frac{r}{2}(\eta_r(\gamma_\nu)-1)\\
\notag&=-\frac{r}{2}\sum\limits_{p=1}^{b} \langle h_p^*,\eta_r(\gamma_\nu)\rangle (\twid{h}_p-1) \pmod{J^2}\\
\label{equation-mod-r-gammas-vs-Kron}&=\frac{r}{2}\sum\limits_{p=1}^{b} \langle h_p^*,\eta_r(\gamma_\nu)\rangle (\twid{h}_p-1) \pmod{J^2}.
\end{align}
The last line follows since in $\Z_r$ for an even $r$, $-\frac{r}{2}=\frac{r}{2}$.

Now, using maps from Mod--$r$ surfaces (i.e. Lemma~\ref{lemma-mod-r-surface-cohomology}) instead of maps from surfaces, we can use the proof from Theorem~\ref{theorem-IntCohom-vs-tors} since \eqref{equation-int-cohom-Fox-just-before-surface} holds for odd $r$, so \eqref{equation-Int-cohom-Fox-after-surface} holds for odd $r$, proving \eqref{equation-Mod-r-CohomFoxAndCupAug} for odd $r$.  For even $r$, \eqref{equation-int-cohom-Fox-just-before-surface} holds with an additional term following from \eqref{equation-Mod-r-aug-Fox-vs-Kron} and \eqref{equation-mod-r-gammas-vs-Kron}.  Specifically, for an even $r$, 
\begin{align*}
\notag\eta_r(\partial r_i &/ \partial x_j)\bmod J^2 =\\
&\sum\limits_{p=1}^b \left ( \sum\limits_\mu \langle h_p^* , \eta_r(\alpha_\mu)
\rangle \langle h_j^* , \eta_r(\beta_\mu) \rangle - \langle h_p^* ,
\eta_r(\beta_\mu) \rangle \langle h_j^* , \eta_r(\alpha_\mu) \rangle\right.\\
&+\frac{r}{2}\left.\sum\limits_\nu\langle h_j^*,\eta_r(\gamma_\nu)\rangle\langle h_p^*,\eta_r(\gamma_\nu)\rangle\right ) (\twid{h}_p - 1).
\end{align*}
This term also occurs in \eqref{equation-Int-cohom-Fox-after-surface} for even $r$ by Lemma~\ref{lemma-mod-r-surface-cohomology}, so \eqref{equation-Mod-r-CohomFoxAndCupAug} holds for even $r$ as well, hence \eqref{equation-Mod-r-CohomFoxAndCup} holds for all $r$.

If we let $a$ be the submatrix of $(\text{pr}\circ\eta)(\partial r_i/\partial x_j)$ consisting of the $b\times b-1$ upper left submatrix,
\[
(h_1-1)\tau(M,e,\omega;r)=|\det(v')|\det(a(1)) = T\det(a(1)) \bmod I^b.
\]
Now by \eqref{equation-Mod-r-CohomFoxAndCup}, computing $\det(a(1))$ is simply computing $q_r(\det(\Theta(1)))$ where $\Theta_{i,j}=\sum\limits_{p=1}^b f_M^r(k_i^*,h_j^*,h_p^*) \twid{h}_p$.
\[
\det(\Theta(1)) = - \twid{h}_1 d(f_M^r,h^*,k^*).
\]
From here, we may follow the proof from Theorem~\ref{theorem-IntCohom-vs-tors}, since $\mu_M^r(h^*,k^*)$ will simply be a sign just as in Theorem~\ref{theorem-IntCohom-vs-tors}, since the linking form of the pseudo-bases is equal to the identity matrix, hence has determinant one.  This follows from equation \eqref{equation-h_i,k_j-linking} which gives for $n+1\leq i\leq b$ and $n\leq j\leq b-1,$ \[\det(h_i\cdot k_j)=\det(\delta_{i,j}) = 1.\] 
\end{proof}
\end{subsection}
%
\end{section}
%
\begin{section}{Integral Massey Products}
\label{section-IntMass}
In this section, we give a generalization of Theorem~\ref{theorem-IntCohom-vs-tors} where we use Massey products rather than the cohomology ring.  The results of this section are similar to results in Chapter~XII Section~2 of \cite{Tur:3d} for closed manifolds.  Note that Massey products are related to Milnor's $\bar\mu$--invariants for links in $S^3$, much like cup products are related to linking numbers (see \cite{Porter-MilMas} and \cite{Tur:MilMas}).
\begin{subsection}{Determinants}
\label{subsection-IntMass-Dets}
First we obtain a new determinant.  Let $R$ be a commutative ring with
$1$, and let $K,L$ be free $R$--modules of rank $n$,$n-1$ respectively,
with $n \geq 2$ and let $S = S(K^*)$, the symmetric algebra on the
dual of $K$, as in Lemma~\ref{lemma-IntCohomDet}.  Let $f\co L\times
K^{m+1} 
\rightarrow R$ be an $R$--map, with $m\geq 1$.  Define $g\co L\times K
\rightarrow S$ by \[g(x,y) = \sum\limits_{i_1,\dots,i_m = 1}^n
f(x,y,a_{i_1},\dots ,a_{i_m})a_{i_1}^*\cdots a_{i_m}^* \in S\] where
$\{a_i\}_{i=1}^n$ is a basis for $K$ and $\{a_i^*\}$ is its dual
basis.  This definition for $g$ looks dependent on the basis chosen, however one can easily show that it is not.

Let $f_0\co L\rightarrow S$ be defined by \[f_0(x) =
\sum\limits_{i_1,\dots,i_{m+1} = 1}^n f(x,a_{i_1},\dots
,a_{i_{m+1}})a_{i_1}^*\cdots a_{i_{m+1}}^* \in S.\]  Again, $f_0$ does
not depend on the chosen basis, by precisely the same argument.  Then
we have the following lemma: 
\begin{lemma}
Suppose $f_0=0$.  Let $a = \{a_i\} , b = \{b_j\}$ be bases of $K,L$
respectively, and let $\theta$ be the $(n-1\times n)$ matrix over $S$
defined by $\theta_{i,j} = g(b_i,a_j)$.  Then there exists a unique $d
= d(f,a,b) \in S^{m(n-1)-1}$ such that 
\begin{equation}
\det(\theta(i)) = (-1)^i a_i^* d.
\label{equation-IntMasseyDet}
\end{equation}
Furthermore, if $a',b'$ are other bases for $K,L$ respectively, then
\begin{equation}
d(f,a',b') = [a'/a][b'/b]d(f,a,b).
\label{equation-IntMasseyCob}
\end{equation}
\label{lemma-IntMasseyDet}
\end{lemma}
\begin{proof}
This is very similar to the proof of Lemma~\ref{lemma-IntCohomDet}.
Let $\beta$ be 
the matrix over $S$ given by $\beta_{i,j} = g(b_i,a_j)a_j^*$.  Then
the sum of the columns of $\beta$ is zero; the $i^\text{th}$ entry in
that sum is $\sum\limits_{j=1}^n \beta_{i,j} = f_0(b_i) = 0$ since our
assumption is $f_0 = 0$.  Now the same argument as given in
Lemma~\ref{lemma-IntCohomDet} to prove \eqref{equation-IntCohomDetDef}
completes the proof of 
\eqref{equation-IntMasseyDet}, and the
argument given to prove \eqref{equation-IntCohomCob} can be used to prove
\eqref{equation-IntMasseyCob}.
\end{proof}
Note that as before, over $\Z$ the determinant is well defined up to
sign, and that one may also sign-refine this determinant to remove the
sign dependence. 

We may also define the condition that $f$ is ``alternate'' in the $K$
variables; let $\overline{f_0}\co L\times K \rightarrow R$ be the $R$--map
given by $\overline{f_0}(x,a) = f(x,\overbrace{a,a,\dots,a}^{m+1
\text{times}})$.  Then $f_0(x) = 0$ for all $x$ clearly implies
$\overline{f_0}(x,a) = 0$ for all $x\in L,a\in K$.  The converse is
also true provided that every polynomial over $R$ which only takes on
zero values has all zero coefficients.  This is true, for example, if
$R$ is infinite with no zero-divisors; in particular for $R=\Z$.  

The rest of the argument is very similar to the argument in \cite{Tur:3d}, section XII.2.  Let $M$ be a $3$--manifold with nonempty
boundary, and for 
$u_1,u_2,\dots , u_k \in H^1(M)$, let $\langle u_1,\dots , u_k\rangle $ denote the
Massey product of $u_1,\dots ,u_k$ as a subset of $H^2(M)$ (note in
general this set may well be empty).  See \cite{Kra:Mas} and
\cite{Fenn:Tech} for definitions and properties of Massey products.
Now assume that $m\geq 1$ is an integer such that
\begin{center}
$(*)_m$:\quad for every $u_1,\dots,u_k\in H^1(M)$ with $k\leq m,
\langle u_1,\dots,u_k\rangle =0$
\end{center}
Here $\langle u_1,\dots,u_k\rangle =0$ means that $\langle u_1,\dots,u_k\rangle$ consists of the
single element $0\in H^2(M)$.  This condition guarantees that for any
$u_1,\dots u_{m+1}\in H^1(M)$, the set $\langle u_1,\dots,u_m\rangle$ consists of a
single element; see \cite{Fenn:Tech} Lemma 6.2.7.  Define a $\Z$--map
$f\co H^1(M,\partial M) \times (H^1(M))^{m+1} \rightarrow \Z$ by 
\[ f(v , u_1 , \dots , u_{m+1}) = (-1)^m\left\langle v \cup \langle u_1,\dots,u_{m+1}\rangle  ,
[M]\right\rangle. \]
The outermost $\langle , \rangle$ is used to denote the evaluation pairing.
\begin{lemma}
$f_0 = 0$, so $f$ has a well-defined determinant (with the sign
refinement as above).
\end{lemma}
For $m=1$, condition $(*)_m$ is void, and in fact the Massey product $\langle u_1,u_2\rangle =-u_1\cup u_2$, so this reduces to Lemma~\ref{lemma-IntCohomDet}.
\begin{proof}
By the argument above, we only need to show that $f$ is alternate.
But this follows from \cite{Kra:Mas} Theorem~15, which gives that for
any element $a\in H^1(M)$, the $m+1$ times Massey product of $a$ with itself,  $\langle\overbrace{a,\dots,a}^{m+1 \text{times}}\rangle,$ lies in
$\Tors(H^2(M))$, hence cupping with an element of $H^1(M,\partial M)$ will give
an element of $\Tors(H^3(M,\partial M))$, which is null.
\end{proof}
We will call this determinant $\Det(f)$, or if we care to introduce the 
sign-refined version with a homology orientation $\omega$, $\Det_{\omega}(f)$.
Since the change of basis formula \eqref{equation-IntMasseyCob} is identical
to the change of basis formula \eqref{equation-IntCohomCob}, the sign
refinement by homology orientation is the same. 
\end{subsection}
%
\begin{subsection}{Relationship to Torsion}
\label{subsection-IntMass-Relns}
\begin{theorem}
Let $M$ be a compact connected oriented $3$--manifold with $\partial M\neq\varnothing,\chi(M)=0$, $n=b_1(M) \geq 2$, and satisfying condition $(*)_m$
for some $m \geq 1$.  Let $e$ be an Euler structure on $M$, let $\omega$ be a homology orientation, and let $q$
be defined as in Section~\ref{section-IntCohom}.  Define the form $f$
as above.  Then 
$\tau(M,e,\omega)\in I^{m(n-1)-1}$ and 
\begin{equation}
\tau(M,e,\omega) \bmod I^{m(n-1)} = q(\Det_\omega(f)) \in
I^{m(n-1)-1}/I^{m(n-1)} .
\label{equation-IntMasseyTors}
\end{equation}
\end{theorem}
\begin{proof}
This proof is very much like the one in Section~\ref{section-IntCohom}.  In place of
\eqref{equation-IntCohomTorsThm}, 
we may use \cite{Tur:MilMas} Theorem~D, which gives the second line of the following string of equalities (all of which are $\bmod J^{m+1}$) {\scriptsize
\begin{align}
\notag \twid{\eta} (\partial r_i/\partial x_j) & = \sum\limits_{i_1,\dots ,
i_m = 1}^n \text{aug} (\twid{\eta}(\partial^{m+1}r_i/\partial
x_{i_1}\dots\partial x_{i_m}\partial x_j))(\twid{h_{i_1}} - 1)\cdots (\twid{h}_{i_m} - 1) \\
\notag & = \sum\limits_{i_1,\dots ,
i_m = 1}^n \left\langle \langle h_{i_1}^* , \dots , h_{i_m}^* , h_j^*\rangle , (-[\Sigma_i])\right\rangle(\twid{h}_{i_1} - 1)\cdots (\twid{h}_{i_m} - 1)\\
\notag & = \sum\limits_{i_1,\dots ,
i_m = 1}^n \left\langle \langle h_{i_1}^* , \dots , h_{i_m}^* , h_j^*\rangle , (-k_i^*\cap [M])\right\rangle (\twid{h}_{i_1} - 1)\cdots (\twid{h}_{i_m} - 1)\\
\notag & = \sum\limits_{i_1,\dots ,
i_m = 1}^n -\left\langle k_i^*\cup \langle h_{i_1}^* , \dots , h_{i_m}^* , h_j^*\rangle , [M]\right\rangle(\twid{h}_{i_1} - 1)\cdots (\twid{h}_{i_m} - 1)\\
\label{equation-reorder-Massey} & = \sum\limits_{i_1,\dots ,
i_m = 1}^n (-1)^{m+1} \left\langle k_i^*\cup \langle h_j^* , h_{i_m}^* , \dots , h_{i_1}^*\rangle , [M]\right\rangle(\twid{h}_{i_1} - 1)\cdots (\twid{h}_{i_m} - 1)\\
\notag & = \sum\limits_{i_1,\dots ,
i_m = 1}^n (-1)^{m+1} \left\langle k_i^*\cup \langle h_j^* , h_{i_1}^* , \dots , h_{i_m}^*\rangle , [M]\right\rangle(\twid{h}_{i_1} - 1)\cdots (\twid{h}_{i_m} - 1)\\
\notag & = \sum\limits_{i_1,\dots ,
i_m = 1}^n -f(k_i^* , h_j^* , h_{i_1}^* , \dots , h_{i_m}^*)(\twid{h}_{i_1} - 1)\cdots (\twid{h}_{i_m} - 1).
\end{align}}
The line marked \eqref{equation-reorder-Massey} follows from \cite{Kra:Mas} Theorem~8, and the next line is by symmetry.  From here, the proof is identical to the proof of Theorem~\ref{theorem-IntCohom-vs-tors} after \eqref{equation-Int-cohom-Fox-after-surface}.
\end{proof}
\end{subsection}
\end{section}
%

\bibliographystyle{amsalpha}
\bibliography{bibList}

\providecommand{\bysame}{\leavevmode\hbox to3em{\hrulefill}\thinspace}
\providecommand{\MR}{\relax\ifhmode\unskip\space\fi MR }
\providecommand{\MRhref}[2]{%
  \href{http://www.ams.org/mathscinet-getitem?mr=#1}{#2}
}
\providecommand{\href}[2]{#2}
\begin{thebibliography}{Tur02}

\bibitem[Fen83]{Fenn:Tech}
Roger~A. Fenn, \emph{Techniques of geometric topology}, London Mathematical
  Society Lecture Note Series, vol.~57, Cambridge University Press, Cambridge,
  1983. \MR{MR787801 (87a:57002)}

\bibitem[Hat02]{Hatcher:AT}
Allen Hatcher, \emph{Algebraic topology}, Cambridge University Press,
  Cambridge, 2002. \MR{MR1867354 (2002k:55001)}

\bibitem[Kra66]{Kra:Mas}
David Kraines, \emph{Massey higher products}, Trans. Amer. Math. Soc.
  \textbf{124} (1966), no.~3, 431--449. \MR{MR0202136 (34 \#2010)}

\bibitem[Nic03]{Nicol-Reid3}
Liviu~I. Nicolaescu, \emph{The {R}eidemeister torsion of 3-manifolds}, de
  Gruyter Studies in Mathematics, vol.~30, Walter de Gruyter \& Co., Berlin,
  2003. \MR{MR1968575 (2004e:57018)}

\bibitem[Por80]{Porter-MilMas}
Richard Porter, \emph{Milnor's {$\bar \mu $}-invariants and {M}assey products},
  Trans. Amer. Math. Soc. \textbf{257} (1980), no.~1, 39--71. \MR{MR549154
  (81a:57021)}

\bibitem[Tur76]{Tur:MilMas}
Vladimir Turaev, \emph{The {M}ilnor invariants and {M}assey products}, Zap.
  Nau\v cn. Sem. Leningrad. Otdel. Mat. Inst. Steklov. (LOMI) \textbf{66}
  (1976), 189--203, 209--210 (Translation in J. Sov. Math. 12 (1979), 128-137),
  Studies in topology, II. \MR{MR0451251 (56 \#9538)}

\bibitem[Tur01]{Tur:Intro}
\bysame, \emph{Introduction to combinatorial torsions}, Lectures in Mathematics
  ETH Z\"urich, Birkh\"auser Verlag, Basel, 2001, Notes taken by Felix Schlenk.
  \MR{MR1809561 (2001m:57042)}

\bibitem[Tur02]{Tur:3d}
\bysame, \emph{Torsions of {$3$}-dimensional manifolds}, Progress in
  Mathematics, vol. 208, Birkh\"auser Verlag, Basel, 2002. \MR{MR1958479
  (2003m:57028)}

\end{thebibliography}
\end{document}